\documentclass[jsco]{academic}
\usepackage{amssymb,amsfonts,amsbsy}
\usepackage{amsmath}
\usepackage{diagrams}
\usepackage{rotating}
\usepackage{diagrams}

\shortauthor{L. Bus\'e and M. Chardin}
\shorttitle{Implicitizing rational hypersufaces}

\title{Implicitizing rational hypersurfaces using approximation
  complexes}

\author[1]{Laurent Bus\'e \thanks{Partially supported by european
    projects GAIA II IST-2001-35512, and
ECG IST-2000-26473.}}
\author[2]{Marc Chardin}

\address[1]{Universit\'e de Nice Sophia-Antipolis,
   Parc Valrose, \\ BP 71,
  06108 Nice Cedex 02, France.\\ lbuse@unice.fr}

\address[2]{Institut de Math\'ematiques, CNRS et Universit\'e Paris 6,\\
4, place Jussieu, F-75252 PARIS CEDEX 05, France.\\chardin@math.jussieu.fr}

\def\KK{{\mathbb{K}}}
\def\PP{{\mathbb{P}}}
\def\ZZ{{\mathbb{Z}}}

\def\Cc{\mathcal{C}}
\def\Sc{\mathcal{S}}
\def\Pc{\mathcal{P}}

\def\f{{\mathbf{f}}}

\def\ux{\underline{x}}
\def\im{{\frak m}}
\def\hom{\hbox{\rm Hom}}
\def\homgr{\hbox{\rm Homgr}}
\def\s{{\sigma}}
\def\om{{\omega}}

\def\sym{\hbox{\rm Sym}}
\def\ann{\hbox{\rm Ann}}

\def\depth{{\mathrm{depth}}}
\def\rank{{\mathrm{rank}}}
\def\indeg{\hbox{\rm indeg}}
\def\Proj{\hbox{\rm Proj}}
\def\proj{\hbox{\rm Proj}}

\def\Zc{{\mathcal{Z}}}
\def\Z{{\mathcal{Z}}}

\begin{document}
\maketitle

\begin{abstract}
  In this paper we describe an algorithm for implicitizing rational 
hypersurfaces in case there exists at most a finite number of base
points. It is based on a technique   exposed in \cite{BuJo02}, 
where implicit equations are obtained as determinants of certain
graded parts of a so-called approximation   complex. We detail and
improve this method by providing an in-depth study of the cohomology
of such a complex. In both particular cases of interest of curve and
surface implicitization we also yield explicit algorithms which only
involves linear algebra routines. 
\end{abstract}

\section{Introduction}

The implicitization problem asks for an implicit equation of a
rational hypersurface given by a
parameterization map $\phi:\PP^{n-1}\rightarrow \PP^{n}$, with $n\geq
2$. This problem received recently a particular interest, especially
in the cases $n=2$ and $n=3$, because it is a key-point in computer
aided geometric design and modeling.   
There are basically three kinds of methods to compute such an implicit
equation of rational curves and surfaces. The first methods use Gr\"obner bases
computations. Even if they always work, they are  
 known to be quite slow in practice and are hence rarely used in
geometric modeling (see e.g. \cite{Hoff89}). The second methods are based on 
resultant matrices. Such methods have the advantage to yield  square
matrices 
whose determinant is an implicit equation. This more compact
formulation of an implicit equation is very useful and well behaved
algorithms are known to work with. However such methods are only
known if there is no base point (see \cite{Jou96}), and in case
base points are isolated and local complete intersection with some
additional technical hypothesis (see \cite{BusPhD}). The third
methods are based on the syzygies of the 
parameterization. They were introduced in \cite{SeCh95}
as the method of ``moving surfaces''. For curve implicitization this
method, also called ``moving lines'', reveals very efficient and
general. This generality is no longer true for surface
implicitization, but the method remains very
efficient. In
\cite{CGZ00} and \cite{DAn01} its validity  was proved
in the absence 
of base point. An extension in the presence of base points is exposed
in \cite{BCD02}, assuming five base points hypothesis, one being that the
base points are isolated and locally complete intersection.

A similar approach involving syzygies was recently explored in
\cite{BuJo02} using more 
systematic tools from algebraic geometry and commutative algebra,
among them the so-called approximation complexes (see
\cite{Vas94}). A new method was given to solve the
hypersurface (and hence curve and surface) implicitization problem in
case base points are isolated and locally 
complete intersection, without any additional hypothesis. 

In this paper we detail and improve the 
results exposed in \cite{BuJo02} for implicitizing rational
hypersurfaces in case the ideal of base points is finite and a
projective local almost complete intersection. We also give some other
properties based on the exactness on these so-called approximation
complexes.  Section 2 and section 3 provide a complete and detailed
description of our results in the particular cases of interest of
respectively curve and surface implicitizations. Algorithms are
completely describe and we present some examples. Section 4 deals with
the general case of hypersurface implicitization, and is much more
technical that both previous one. In particular all the proofs are
given in this section. 
\medskip

Hereafter $\KK$ denotes any field.

\section{Implicitization of rational parametric curves}
In this section we present a method to implicitize any parameterized 
plane curve. Let $f_0,f_1,f_2$ be three
homogeneous polynomials in $\KK[s,t]$ of the same degree $d\geq 1$, and
consider the rational map
\begin{eqnarray*}
\phi :\ \ \ \ \ \PP^1_\KK & \rightarrow & \PP^2_\KK \\
   (s:t) & \mapsto & (f_0(s,t):f_1(s,t):f_2(s,t)).
\end{eqnarray*}
We denote by $x,y,z$ the homogeneous coordinates of $\PP^2$.  The
algebraic closure of the image of $\phi$ is a curve in $\PP^2$ if and
only if $\phi$ is generically finite onto its image, that is
$\delta:=\gcd(f_0,f_1,f_2)$ is not of degree $d$. We assume this minimal
hypothesis and denote by $\Cc$ the scheme-theoretic closed image
of $\phi$. It is known that $\Cc$ is an irreducible and reduced curve of degree $(d-\deg(\delta))/\beta$, 
where $\beta$ denotes the degree
of $\phi$ onto its image (that is the number of points in a general
fiber). Thus any implicit equation of $\Cc$ is of the form $P(x,y,z)$ 
where $P$ is an irreducible homogeneous polynomial of degree 
$(d-\delta)/\beta$. 

In what follows we describe an algorithm to compute explicitly an
implicit equation $P$ of $\Cc$ (in fact we will compute $P^\beta$) 
without any other hypothesis on the 
polynomials $f_0,f_1,f_2$. In the case $\deg(\delta)=0$ we recover
the well-known method of moving lines (see \cite{SeCh95,CSC98}).

\subsection{The method}

Let us denote by $A$ the polynomial ring $\KK[s,t]$. We consider it
with its natural graduation obtained by setting $\deg(s)=\deg(t)=1$.
From the polynomials $f_0,f_1,f_2$ of the parameterization $\phi$, we
can build the following well-known graded Koszul complex (the notation
$[-]$ stands for the degree shift in $A$) 
\begin{equation}\label{koszul_A}
0 \rightarrow A[-3d] \xrightarrow{d_3} A[-2d]^3 \xrightarrow{d_2}
A[-d]^3 \xrightarrow{d_1} A,
\end{equation}
where the differentials are given by
$$d_3=\begin{pmatrix} f_2 \\ -f_1 \\ f_0 \end{pmatrix},
d_2=\begin{pmatrix}
  f_1 & f_2 &  0 \\
  -f_0 & 0 & f_2 \\
  0 & -f_0 & -f_1 \end{pmatrix}, d_1=\begin{pmatrix} f_0 & f_1 & f_2
\end{pmatrix}.$$
In what follows we will not consider exactly this complex, but the
complex obtained by tensorizing it by $A[x,y,z]$ over $A$. This
complex, that we denote  $(K_\bullet(f_0,f_1,f_2),u_\bullet)$, is of
the form
$$
0 \rightarrow A[x,y,z][-3d] \xrightarrow{u_3} A[x,y,z][-2d]^3
\xrightarrow{u_2} A[x,y,z][-d]^3 \xrightarrow{u_1} A[x,y,z],$$
where
the matrices of the differentials $d_i$ and $u_i$ are the same, for all 
$i=1,2,3$.  Note that the ring $A[x,y,z]$ is naturally
bi-graded, having a graduation coming from $A=\KK[s,t]$, and another
one coming from $\KK[x,y,z]$ with $\deg(x)=\deg(y)=\deg(z)=1$; we  
hereafter adopt the notation $(-)$ for the degree shift in 
$\KK[x,y,z]$. 

We form another bi-graded Koszul complex on $A[x,y,z]$, the one
associated to the sequence $(x,y,z)$. We denote it by
$(K_\bullet(x,y,z),v_\bullet)$, it is of the form
$$
0 \rightarrow A[x,y,z](-3) \xrightarrow{v_3} A[x,y,z](-2)^3
\xrightarrow{v_2} A[x,y,z](-1)^3 \xrightarrow{v_1} A[x,y,z],$$
and the matrices of its differentials are obtained from the matrices
of the differentials of \eqref{koszul_A} by replacing $f_0$ by $x$,
$f_1$ by $y$ and $f_2$ by $z$. 
Observe that since $(x,y,z)$ is a regular sequence in $A[x,y,z]$, the
previous complex $K_\bullet(x,y,z)$ is acyclic, that is to say all its
homology groups $H_{i}(K_\bullet(x,y,z))$ vanish for $i>0$.

We can now construct a new bi-graded complex of $A[x,y,z]$-modules, denoted $\Zc_\bullet$, from both
Koszul complexes $(K_\bullet(f_0,f_1,f_2),u_\bullet)$ and
$(K_\bullet(x,y,z),v_\bullet)$ (observe that these complexes differ
only by their differentials). Define $Z_i:=\ker(d_i)$
for all $i=0,\ldots,3$ (with $d_0:A\rightarrow 0$), and set 
$\Zc_i:=Z_i[id]\otimes_A A[x,y,z]$ for $i=0,\ldots,3$, which are bi-graded $A[x,y,z]$-modules. The map $v_1$ induced obviously the bi-graded map, that we denote also by $v_1$,   
\begin{eqnarray*}
  \Zc_1(-1) & \xrightarrow{v_1} & \Zc_0=A[x,y,z] \\
  (g_1,g_2,g_3) & \mapsto & g_1x+g_2y+g_3z. 
\end{eqnarray*}
Using the differential
$v_2$ of  $K_\bullet(x,y,z)$ we can map $\Zc_2$ to
$A[x,y,z](-d)^3$, but since $u_1\circ v_2 + v_1 \circ u_2=0$ (which
follows from a straightforward computation), we have $v_2(\Zc_2)
\subset \Zc_1$. And in the same way we can map $\Zc_3$ to $\Zc_2$ with the differential $v_3$, since $u_2\circ v_3 + v_2\circ u_3=0$. Thus we obtain the following 
bi-graded complex (it is a complex since
$(K_\bullet(x,y,z),v_\bullet)$ is)~:
$$(\Zc_\bullet,v_\bullet)~: 0 \rightarrow \Zc_3(-3)
\xrightarrow{v_3} \Zc_2(-2) \xrightarrow{v_2} \Zc_1(-1)
\xrightarrow{v_1} \Zc_0=A[x,y,z].$$
This complex is known as the \emph{
approximation complex of cycles} associated to the polynomials $f_0,f_1,f_2$ in
$\KK[s,t]$. It was originally introduced in \cite{SiVa81} for studying
Rees algebras through symmetric algebras (see also \cite{Vas94}).

\begin{remark}
  Point out that in the language of the moving lines method, 
the image of a
  triple $(g_1,g_2,g_3) \in {\Zc_1}_{[\nu](0)}$ by $v_1$ is nothing 
but a moving 
  line of degree $\nu$ following the surface parameterized by
  the polynomials $f_0,f_1,f_2$ (see e.g. \cite{Cox01}).
\end{remark}

Since we have supposed that $f_0,f_1,f_2$ are all non-zero we have
$Z_3=0$, and consequently $\Zc_3=0$, and the following theorem (recall $
\delta:=\gcd(f_0,f_1,f_2)$)~: 

\begin{theorem}\label{theorem curve}  The determinant of the
  graded complex of free $\KK[x,y,z]$-modules
  $$0 \rightarrow {\Zc_2}_{[d-1]}(-2) \xrightarrow{v_2}
  {\Zc_1}_{[d-1]}(-1) \xrightarrow{v_1}
  {\Zc_0}_{[d-1]}=A[x,y,z]_{[d-1]},$$
  is $P(x,y,z)^\beta$, where $P$ is an implicit
  equation of the curve $\Cc$.
  
  Moreover, if $\deg(\delta)=0$ then ${\Zc_2}_{[d-1]}=0$; 
  thus the $d\times d$ determinant of the 
  map ${\Zc_1}_{[d-1]}(-1) \xrightarrow{v_1}
  {\Zc_0}_{[d-1]}$ equals $P(x,y,z)^\beta$.
\end{theorem}
\begin{proof}
  See section 5 in \cite{BuJo02} for a proof. We only mention how
  the last statement follows from the description of the approximation
  complex $\Zc_\bullet$. If $\deg(\delta)=0$ then
  $\depth_{(s,t)}(f_0,f_1,f_2)=2$, which means that $f_0,f_1,f_2$ have
  no base points in $\PP^1$. This implies that not only the third
  homology group of the Koszul complex \eqref{koszul_A} vanishes, but
  also the second. It follows that $Z_2\simeq A[-3d]$, and hence 
$\Zc_2\simeq A[x,y,z][-d]$.
\end{proof}

The second statement of this theorem gives exactly the matrix
constructed by the method of moving lines. The first statement show
that this  method can be extended even if we do not assume
$\deg(\delta)=0$, $P^\beta$ being obtained as the quotient of two determinants of respective size $d$ and $\delta$.  
In \cite{BuJo02} it is in fact proved that for any integer $\nu\geq
d-1$ the determinant of the complex 
$(\Zc_\bullet)_{[\nu]}$ equals $P^\beta$. In case
$\deg(\delta)=0$, this and the graded isomorphism $\Zc_2\simeq
A[x,y,z][-d]$ explain clearly, in our point of view, 
why the method of moving lines works so well with and \emph{only} with moving
lines of degree $d-1$. 

\subsection{The algorithm}
We here convert theorem \ref{theorem curve} into an explicit algorithm where each step reduces to well-known and efficient linear algebra routines.

\medskip

\noindent {\sc Algorithm} (for implicitizing rational parametric
curves):

\noindent {\bf Input~:} Three homogeneous polynomials
$f_0(s,t),f_1(s,t),f_2(s,t)$ of the same degree $d\geq 1$ such that
$\delta=\gcd(f_0,f_1,f_2)$ is not of degree $d$.

\noindent {\bf Output~:} Either~:
\begin{itemize}
\item a square matrix $\Delta_1$ such that $\det(\Delta_1)$ equals $P^\beta$, in case $\deg(\delta)=0$,
\item two square matrices $\Delta_1$ and $\Delta_2$, respectively of
  size $d$ and $\deg(\delta)$, such that
  $\frac{\det(\Delta_1)}{\det(\Delta_2)}$ equals $P^\beta$.
\end{itemize}

{\sl \begin{enumerate}
  \item Compute the matrix ${\tt F_1}$ of the first map of
    \eqref{koszul_A}~: $A^3_{d-1}\xrightarrow{d_1} A_{2d-1}$. Its
    entries are either 0 or a coefficient of $f_0,f_1$ or $f_2$, and
    it is of size $2d\times 3d$.
  \item Compute a kernel matrix ${\tt K_1}$ of the transpose of ${\tt F_1}$. It has
    $3d$ columns and $\rank({\Zc_1}_{[d-1]})$ lines.
  \item Construct the matrix ${\tt Z_1}$ by
    ${\tt Z_1}(i,j)=x{\tt K_1}(j,i)+y{\tt K_1}(j,i+d),z{\tt K_1}(j,i+2d)$, with $i=1,\ldots,d$
    and $j=1,\ldots,\rank({\Zc_1}_{[d-1]})$. It is a matrix of the
    map ${\Zc_1}_{[d-1]}(-1) \xrightarrow{v_1} {\Zc_0}_{[d-1]}$.
  \item \textbf{If} ${\tt Z_1}$ is square \textbf{then} set
    $\Delta_1:={\tt Z_1}$ \textbf{else} 
    \begin{enumerate}
  \item Compute a list $L_1$ of $d$ integers indexing $d$
    independent columns in ${\tt Z_1}$. Let $\Delta_1$ be the $d\times d$
    submatrix of ${\tt Z_1}$ obtained by removing columns not in $L_1$.
  \item Compute the matrix ${\tt F_2}$ of the second map of
    \eqref{koszul_A}~: $A^3_{d-1}\xrightarrow{d_2} A^3_{2d-1}$, and
    a kernel matrix ${\tt K_2}$ of its transpose. The matrix ${\tt K_2}$ has $3d$
    columns and $\rank({\Zc_2}_{[d-1]})$ lines.
  \item Construct the matrix ${\tt Z'_2}$ by, for all
    $j=1,\ldots,\rank({\Zc_2}_{[d-1]})$,
   \begin{eqnarray*}
     i=1,\ldots,d & : & {\tt Z'_2}(i,j)=y{\tt K_2}(j,i)+z{\tt K_2}(j,i+d), \\
     i=d+1,\ldots,2d & : & {\tt Z'_2}(i,j)=-x{\tt K_2}(j,i-d)+z{\tt K_2}(j,i+d), \\
     i=2d+1,\ldots,3d & : & {\tt Z'_2}(i,j)=-x{\tt K_2}(j,i-d)-y{\tt K_2}(j,i).
   \end{eqnarray*}
 \item Construct the $\rank({\Zc_1}_{[d-1]}) \times
   \rank({\Zc_2}_{[d-1]})$ matrix ${\tt Z_2}$ whose $j^{\mathrm{th}}$ column
   ${\tt Z_2}(\bullet,j)$ is the solution of the linear system
   ${}^t{\tt Z_2}(\bullet,j).{\tt K_1}={\tt Z'_2}(\bullet,j)$. It is a matrix of the
   map ${\Zc_2}_{[d-1]}(-1) \xrightarrow{v_2} {\Zc_1}_{[d-1]}$.
 \item Define $\Delta_2$ to be the square submatrix of ${\tt Z_2}$
   obtained by removing the lines indexed by $L_1$.
\end{enumerate}
\textbf{endif}
 \end{enumerate}}

\begin{remark} In order to keep this algorithm easily understandable
  we kept the symbolic variables $t$ and $z$, but they are of course
  useless and should be specialized to 1.
\end{remark}

\subsection{An example}

As we have already said, the previous algorithm is exactly the well-known 
method of moving 
lines in case $\deg(\delta)=0$. The only thing new is its
extension to the case where $\deg(\delta)>0$. We illustrate it with the
following (very simple) example.

Let $f_0(s,t)=s^2$, $f_1(s,t)=st$ and $f_2(s,t)=t^2$. Applying our
algorithm we found that the matrix ${\tt Z_1}$ (step 3) is square~:
$${\tt Z_1}=\begin{pmatrix}
  -y & -z \\
  x & y \end{pmatrix}.$$
We deduce that an implicit equation is given by $xz-y^2$. Now let us
multiply artificially each $f_i$ by $s$. We obtain new entries for our 
algorithm
which are~: $f_0(s,t)=s^3$, $f_1(s,t)=s^2t$ and
$f_2(s,t)=st^2$. In this case the matrix ${\tt Z_1}$ is not square, it is
the $3\times 4$ matrix~:
$${\tt Z_1}=\begin{pmatrix}
  -y & -z & -z & 0 \\
  x & 0 & y & -z \\
 0 & x & 0 & y \end{pmatrix}.$$
The matrix $\Delta_1$ can be chosen to be the submatrix of ${\tt Z_1}$ given by $\begin{pmatrix}
  -y & -z & -z  \\
  x & 0 & y \\
 0 & x & 0  \end{pmatrix},$ with determinant $-x^2z + xy^2$. 
Continuing the algorithm  we obtain the matrices ${\tt Z_2}=\begin{pmatrix}
  -z \\ y \\ 0 \\ -x \end{pmatrix}$, and $\Delta_2=\begin{pmatrix} -x
\end{pmatrix}$. It follows that an implicit equation is $-xz+y^2$.

\section{Implicitization of rational parametric surfaces}\label{surf}
We arrive to the much more intricate problem of surface implicitization. Let
$f_0,f_1,f_2,f_3$ be four homogeneous polynomials in $\KK[s,t,u]$ of
the same degree $d\geq 1$, and consider the rational map
\begin{eqnarray*}
\phi :\ \ \ \ \ \PP^2_\KK & \rightarrow & \PP^3_\KK \\
   (s:t:u) & \mapsto & (f_0(s,t,u):f_1(s,t,u):f_2(s,t,u):f_3(s,t,u)).
\end{eqnarray*}
We denote by $x,y,z,w$ the homogeneous coordinates of $\PP^2$.  The
closure of the image of $\phi$ is a surface if and only if the map
$\phi$ is generically finite onto its image, which we assume hereafter. Thus
let $\Sc$ denote the surface in $\PP^3$ obtained as the
 closed image of $\phi$. Our problem is to compute
an implicit equation of $\Sc$. It is more difficult than curve
implicitization because of the presence of base points in codimension
2 arising from the parameterization, base points in codimension 1
being easily removed by substituting each $f_i$, $i=0,\ldots,3$, by 
$f_i/\delta$, where $\delta:=
\gcd(f_0,f_1,f_2,f_3)$, if necessary. Consequently, we assume that the
ideal $I=(f_0,f_1,f_2,f_3)$ in $\KK[s,t,u]$ is at least of codimension 2, that
is defines $n$ 
isolated points $p_1,\ldots,p_n$ in $\PP^2$, and we have (see section
2 and theorem 2.5 in \cite{BuJo02} for a proof and more details)~:
\begin{theorem} Let $e_{p_i}$ denote the algebraic 
  multiplicity of the point $p_i$ for all $i=1,\ldots,n$, and $\beta$
  denote the degree of $\phi$ 
  onto its image. Then
  $$d^2-\sum_{i=1}^n e_{p_i}=\left\{\begin{array}{cl} \beta\deg(\Sc) &
      \ \ \mathrm{if} \ \phi
      \ \mathrm{is \ generically \ finite,}\\
      0 & \ \ \mathrm{if} \ \phi \ \mathrm{is \ not \ generically \ 
        finite.}
        \end{array}\right.$$
\end{theorem}
This theorem gives us the degree of the closed image of $\phi$ if it
is a surface (that we have denoted by $\Sc$), and also that this image is not
a surface if and only if $d^2=\sum_{i=1}^n e_{p_i}$.  We will call an
\emph{implicit equation} of $\Sc$ an equation of its associated
divisor, which is an irreducible and homogeneous polynomial 
$P(x,y,z,w)$ of total degree $(d^2-\sum_{i=1}^n e_{p_i})/\beta$.

In the following subsections we mainly present a method based on the
approximation complexes to compute explicitly an implicit equation of $\Sc$ (in fact we will compute $P^\beta$) in
case the ideal $I$ is a local complete intersection of codimension at least 2, and we also give some other related results. Then we expose an explicit description of the algorithm it yields and illustrate it with some examples.

\subsection{The method} 
We denote by $A$ the polynomial ring
$\KK[s,t,u]$ which is naturally graded by $\deg(s)=\deg(t)=\deg(u)=1$. 
Let us form the Koszul complex of $f_0,f_1,f_2,f_3$ in $A$~:
\begin{equation}\label{Koszul_A}
0 \rightarrow A[-4d] \xrightarrow{d_4} A[-3d]^4 \xrightarrow{d_3}
A[-2d]^6 \xrightarrow{d_2} 
A[-d]^4 \xrightarrow{d_1} A,
\end{equation}
where the differentials are given by
$$d_4=\begin{pmatrix} -f_3 \\ f_2 \\ -f_1 \\f_0 \end{pmatrix}, \ \ 
d_3=\begin{pmatrix}
  f_2 & f_3 & 0 & 0  \\
  -f_1 & 0 & f_3 & 0 \\
  f_0 & 0 & 0 & f_3 \\
  0 & -f_1 & -f_2 & 0 \\
  0 & f_0 & 0 & -f_2 \\
  0 & 0 & f_0 & f_1 \end{pmatrix},$$
$$d_2=\begin{pmatrix}
  -f_1 & -f_2 &  0 & -f_3 & 0 & 0 \\
  f_0 & 0 & -f_2 & 0 & -f_3 & 0 \\
  0 & f_0 & f_1 & 0 & 0 & -f_3 \\
  0 & 0 & 0 & f_0 & f_1 & f_2 \end{pmatrix}, \ \ 
d_1=\begin{pmatrix} f_0 & f_1 & f_2 & f_3 \end{pmatrix}.$$
As for curve implicitization, we denote by $(K_\bullet(f_0,f_1,f_2,f_3),u_\bullet)$ this
Koszul complex tensorized by $A[\ux]:=A[x,y,z,w]$ over $A$,
which is of the form~:
$$0 \rightarrow A[\ux][-4d] \xrightarrow{u_4} A[\ux][-3d]^4
\xrightarrow{u_3} A[\ux][-2d]^6 
\xrightarrow{u_2} A[\ux][-d]^4 \xrightarrow{u_1} A[\ux],$$
where the matrices of the differentials $d_i$ and $u_i$ are the same,
$i=1,2,3,4$ (here again we set  $\deg(x)=\deg(y)=\deg(z)=\deg(w)=1$).
We also consider the bi-graded Koszul complex on $A[\ux]$ 
associated to the sequence $(x,y,z,w)$, and denote it 
$(K_\bullet(x,y,z,w),v_\bullet)$~:
$$0 \rightarrow A[\ux](-4) \xrightarrow{v_4} A[\ux](-3)^4
\xrightarrow{v_3} A[\ux](-2)^6 
\xrightarrow{v_2} A[\ux](-1)^4 \xrightarrow{v_1} A[\ux].$$
The matrices of its differentials are obtained from the matrices of the
differentials of \eqref{Koszul_A} by replacing $f_0$ by $x$, $f_1$ by
$y$, $f_2$ by $z$ and $f_3$ by $w$. Note that since $(x,y,z,w)$
is a regular sequence in $A[\ux]$, the 
complex $K_\bullet(x,y,z,w)$ is acyclic.
From both Koszul complexes $(K_\bullet(f_0,f_1,f_2,f_3),u_\bullet)$ and
$(K_\bullet(x,y,z,w),v_\bullet)$ we can build, as we did for curves, the  approximation complex $\Zc_\bullet$. We define $Z_i:=\ker(d_i)$ and $\Zc_i:=Z_i[id]\otimes_A A[\ux]$ for all
$i=0,1,2,3,4$ (where $d_0:A\rightarrow 0$), they are naturally bi-graded $A[\ux]$-modules. Since for all
$i=1,2,3$ we have $u_i\circ v_{i+1} + v_i \circ u_{i+1}=0$,  we obtain
the  bi-graded complex ~:
$$(\Zc_\bullet,v_\bullet)~: 0 \rightarrow \Zc_4(-4)
\xrightarrow{v_4}  \Zc_3(-3)
\xrightarrow{v_3} \Zc_2(-2) \xrightarrow{v_2} \Zc_1(-1)
\xrightarrow{v_1} \Zc_0=A[\ux],$$
where, since we have supposed $d\geq 1$, $\Zc_4=0$.

\begin{remark} In the language of the moving surfaces method
  (see \cite{SeCh95,Cox01}), an element 
  $(g_1,g_2,g_3,g_4) \in {\Zc_1}_{[\nu](0)}$ is nothing but a moving
  hyperplane of degree $\nu$ following the surface $\Sc$.
\end{remark}

To state the main result of this section, we need some notations. If $p$ is an isolated base point defined by the ideal $I$, we denote by $d_p$ its geometric multiplicity (also called its degree); note that we have already denoted by $e_p$ its algebraic multiplicity. Recall that if $M$ is a $\ZZ$-graded $R$-module, where $R$ is a $\ZZ$-graded ring, its initial degree is defined as $\indeg(M)=\min\{\nu \in \ZZ : M_\nu\neq 0 \}$.

\begin{theorem}\label{theorem surface} Suppose that the ideal
  $I=(f_0,f_1,f_2,f_3)\subset A$ is of codimension at least 2 and $\KK$ is infinite. Let $\mathcal{P}:=\Proj(A/I)$ and denote by $I_\mathcal{P}$ the saturated ideal of $I$ w.r.t. the maximal ideal $\mathfrak{m}=(s,t,u)$ of $A$. Then,
\begin{itemize}
 \item The complex $\Zc_\bullet$ is acyclic if and only if $\mathcal{P}$ is locally generated by (at most) $3$ elements. \\
 \item Assume that $\mathcal{P}$ is locally generated by  $3$ elements, then for all integer $$\nu \geq \nu_0:=2(d-1)-\indeg(I_\mathcal{P})$$ 
the determinant of the 
  graded complex of free $\KK[\ux]$-modules
  $$0 \rightarrow {\Zc_3}_{[\nu]}(-3) \xrightarrow{v_3} 
  {\Zc_2}_{[\nu]}(-2) \xrightarrow{v_2}
  {\Zc_1}_{[\nu]}(-1) \xrightarrow{v_1}
  {\Zc_0}_{[\nu]}=A_{[\nu]}[\ux]$$
 is an homogeneous element of $\KK[\ux]$ of degree $d^2-\sum_{p\in
   \mathcal{P}}d_p$,  and is a multiple of $P^\beta$ independent of $\nu$, where $P$ is an implicit equation of
 $\Sc$. It is exactly $P^\beta$ if
 and only if $I$ is locally a complete intersection. Moreover, for all
 $\nu\in \ZZ$, ${\Zc_3}_{[\nu]}$ is always a free $\KK[\ux]$-module of
 rank $\max(\binom{\nu-d+2}{2},0)$; in particular ${\Zc_3}_{[\nu]}=0$
 if and only if $\nu \leq d-1$.  
\end{itemize}
\end{theorem}

\begin{proof}
  See theorem \ref{TH} for the complete proof. Note that an argument similar
    to the one given in the proof of theorem \ref{theorem
    curve} shows easily that ${\Zc_3}\simeq A[\ux][-d]$, and hence the
    last statement of the second point of the theorem. 
\end{proof}

By standard properties of determinants of complexes (see e.g. appendix
A in \cite{GKZ94}) we deduce the 
\begin{corollary}\label{Mf} 
Suppose that $I=(f_0,f_1,f_2,f_3)\subset A$ is of codimension at least 2, $\KK$ is infinite, and $\Pc=\Proj(A/I)$ is locally generated by 3 elements. Then 
for all $\nu\geq \nu_0:=2(d-1)-\indeg(I_\Pc)$ any non-zero minor of (maximal) size $(\nu+2)(\nu+1)/2$ of
  the surjective matrix 
  \begin{eqnarray*}
    {\Zc_1}_{[\nu]}(-1) & \xrightarrow{v_1} & A_{[\nu]}[\ux] \\
    (g_1,g_2,g_3,g_4) & \mapsto & xg_1+yg_2+zg_3+wg_4
   \end{eqnarray*}
  is a non-zero multiple of $P^\beta$. Moreover, if $\Pc$ is locally a complete intersection, then 
  the gcd of all these minors equals $P^\beta$.
\end{corollary}

Before going through an algorithmic version of this theorem we make few
remarks on the integer $\nu_0$. First note that $\nu_0$ depends geometrically on the ideal $I$,
since $\indeg(I_\Pc)$ is the smallest degree of a hypersurface in $\PP^2$ containing the
closed subscheme defined by $I$.  Let us observe how the initial degree
of $I_\Pc$ behaves. If $I$ has no base points, that is $I_\Pc=A$, then
$\indeg(I_\Pc)=0$. If there exists base points then this initial
degree is always greater or equal to 1 since $I_\Pc$ is generated in
degree at least 1, and is always bounded by $d$ since the $f_i$'s are in $I_\Pc$. Also if $I$ is saturated
then its initial degree is exactly $d$. We deduce that 
\begin{itemize}
\item $\nu_0=2d-2$ if $I$ has no base points,
\item $d-2\leq \nu_0\leq 2d-3$ if $I$ has base points,
\item $\nu_0=d-2$ if $I$ is saturated (note that in this case we know that ${\Zc_3}_{[\nu_0]}=0$, and hence $\det({\Zc_\bullet}_{[\nu_0]})$ is always obtained as a single determinant or as a quotient of two determinants).
\end{itemize} 
This shows in particular that the presence of base points simplify the complexity of the computation of the implicit surface.  Finally recall that the explicit  bound (and not the theoretical one) given in \cite{BuJo02} is only $2d-2$, if there exists base points or not.  

\subsection{The algorithm}

We now  develop the algorithm suggested by theorem
\ref{theorem surface}, and then discuss some computational aspects.

\medskip

\noindent {\sc Algorithm} (for implicitizing rational parametric
surfaces with local almost complete intersection isolated base points, possibly empty):

\noindent {\bf Input~:} Four homogeneous polynomials
$f_0,f_1,f_2,f_3$ in $A$ of the same degree $d\geq 1$ such that
the ideal $(f_0,f_1,f_2,f_3)\subset A$ is  locally generated by 3 elements 
outside $V(s,t,u)$ and at least of codimension 2. An integer $\nu$ with default value
$\nu:=2d-2$. 

\noindent {\bf Output~:} Either~:
\begin{itemize}
\item a square matrix $\Delta_1$ such that $\det(\Delta_1)$ is the determinant of ${(\Zc_\bullet)}_{[\nu]}$,
\item two square matrices $\Delta_1$ and $\Delta_2$ such that
  $\frac{\det(\Delta_1)}{\det(\Delta_2)}$ is the determinant of ${(\Zc_\bullet)}_{[\nu]}$,
\item three square matrices $\Delta_1,\Delta_2$ and
$\Delta_3$ such that
$\frac{\det(\Delta_1)\det(\Delta_3)}{\det(\Delta_2)}$ is the determinant of ${(\Zc_\bullet)}_{[\nu]}$.
\end{itemize}

{\sl \begin{enumerate}
  \item Compute the matrix ${\tt F_1}$ of the first map of
    \eqref{Koszul_A}~: $A^4_{\nu}\xrightarrow{d_1} A_{\nu+d}$, and a 
    kernel ${\tt K_1}$ of its transpose which has $\rank({\Zc_1}_{[\nu]})$
    lines. 
  \item Set $m:=\frac{(\nu+2)(\nu+1)}{2}$. Construct the matrix ${\tt Z_1}$
    defined by~:
    {\small $${\tt Z_1}(i,j)=x{\tt K_1}(j,i)+y{\tt K_1}(j,i+m),z{\tt K_1}(j,i+2m)+w{\tt K_1}(j,i+3m),$$} 
    with $i=1,\ldots,m$   
    and $j=1,\ldots,\rank({\Zc_1}_{[\nu]})$.
    It is the matrix of the
    map ${\Zc_1}_{[\nu]}(-1) \xrightarrow{v_1} {\Zc_0}_{[\nu]}$.
  \item \textbf{If} ${\tt Z_1}$ is square \textbf{then} set
    $\Delta_1:={\tt Z_1}$
    \textbf{else}
    \begin{enumerate}
  \item Compute a list $L_1$ of integers indexing 
    independent columns in ${\tt Z_1}$. $L_1$ consists in $m$ integers. Let
    $\Delta_1$ be the $m\times m$ 
    submatrix of ${\tt Z_1}$ obtained by removing columns not in $L_1$.
  \item Compute the matrix ${\tt F_2}$ of the second map of
    \eqref{Koszul_A}~: $A^6_{\nu}\xrightarrow{d_2} A^3_{\nu+d}$, and
    a kernel ${\tt K_2}$ of its transpose which has
    $\rank({\Zc_2}_{[\nu]})$ lines. 
  \item Construct the matrix ${\tt Z'_2}$ defined by, for all
    $j=1,\ldots,\rank({\Zc_2}_{[\nu]})$,
   {\small $$\begin{array}{l}
     {\tt Z'_2}(i,j)=-y{\tt K_2}(j,i)-z{\tt K_2}(j,i+m)-w{\tt K_2}(j,i+3m), i=1,\ldots,m,\\
     {\tt Z'_2}(i,j)=x{\tt K_2}(j,i-m)-z{\tt K_2}(j,i+m)-w{\tt K_2}(j,i+3m), i=m+1,\ldots,2m,\\
     {\tt Z'_2}(i,j)=x{\tt K_2}(j,i-m)+y{\tt K_2}(j,i)-w{\tt K_2}(j,i+3m), i=2m+1,\ldots,3m,\\
     {\tt Z'_2}(i,j)=x{\tt K_2}(j,i)+y{\tt K_2}(j,i+m)+z{\tt K_2}(j,i+2m), i=3m+1,\ldots,4m.
   \end{array}$$}
 \item Construct the $\rank({\Zc_1}_{[\nu]}) \times
   \rank({\Zc_2}_{[\nu]})$ matrix ${ \tt Z_2}$ whose $j^{\mathrm{th}}$ column
   ${\tt Z_2}(\bullet,j)$ is the solution of the linear system
   ${}^t{\tt Z_2}(\bullet,j).{\tt K_1}={\tt Z'_2}(\bullet,j)$. It is the matrix of the
   map ${\Zc_2}_{[\nu]}(-1) \xrightarrow{v_2} {\Zc_1}_{[\nu]}$.
 \item Define $\Delta'_2$ to be the submatrix of ${\tt Z_2}$
   obtained by removing the lines indexed by $L_1$. \\
   \textbf{If} $\Delta'_2$ is square \textbf{then} set
    $\Delta_2:=\Delta'_2$
    \textbf{else}
    \begin{enumerate}  
    \item Compute a list $L_2$ of integers indexing independent columns in
   $\Delta'_2$. Define $\Delta_2$ to be the square submatrix of
   $\Delta'_2$ obtained by removing columns not in $\L_2$.
    \item Construct the matrix ${\tt F_3}$ of the third map of
    \eqref{Koszul_A}~: $A^4_{\nu}\xrightarrow{d_3} A^6_{\nu+d}$, and
    the kernel ${\tt K_3}$ of its transpose which has
    $\rank({\Zc_3}_{[\nu]})$ lines.  
    \item Construct the matrix ${\tt Z'_3}$ defined by, for all
    $j=1,\ldots,\rank({\Zc_3}_{[\nu]})$,
    {\small $$\begin{array}{l}
     {\tt Z'_3}(i,j)=z{\tt K_3}(j,i)+w{\tt K_3}(j,i+m), i=1,\ldots,m,\\
     {\tt Z'_3}(i,j)=-y{\tt K_3}(j,i-m)+w{\tt K_3}(j,i+m), i=m+1,\ldots,2m,\\
     {\tt Z'_3}(i,j)=x{\tt K_3}(j,i-2m)+w{\tt K_3}(j,i+m), i=2m+1,\ldots,3m,\\
     {\tt Z'_3}(i,j)=-y{\tt K_3}(j,i-2m)-z{\tt K_3}(j,i-m), i=3m+1,\ldots,4m,\\
     {\tt Z'_3}(i,j)=x{\tt K_3}(j,i-3m)-z{\tt K_3}(j,i-m), i=4m+1,\ldots,5m,\\
     {\tt Z'_3}(i,j)=x{\tt K_3}(j,i-3m)+y{\tt K_3}(j,i-2m), i=5m+1,\ldots,6m.
   \end{array}$$}
   \item Construct the $\rank({\Zc_2}_{[\nu]}) \times
   \rank({\Zc_3}_{[\nu]})$ matrix ${\tt Z_3}$ whose $j^{\mathrm{th}}$ column
   ${\tt Z_3}(\bullet,j)$ is the solution of the linear system
   ${}^t{\tt Z_3}(\bullet,j).{\tt K_2}={\tt Z'_3}(\bullet,j)$. It is the matrix of the
   map ${\Zc_3}_{[\nu]}(-1) \xrightarrow{v_3} {\Zc_2}_{[\nu]}$.
  Define $\Delta_3$ to be the square submatrix of ${\tt Z_3}$
   obtained by removing the lines indexed by $L_2$. $\Delta_3$ is of
   size $\frac{(\nu-d+2)(\nu-d+1)}{2}$. 
\end{enumerate}   
\textbf{endif}
 \end{enumerate}
 \textbf{endif.}
 \end{enumerate}}

In the language of the moving surfaces method, the matrix ${\tt Z_1}$ obtained at step 2 
gather all the moving \emph{hyperplanes} of degree $\nu$ following
the surface $\Sc$. Assume hereafter that $\Pc$ is locally a complete intersection. From a computational point of view, corollary
\ref{Mf} implies that this matrix can be taken as a
\emph{representation} of the surface $\Sc$, replacing an expanded
implicit equation (even if it is generally non-square). For
instance, to test if a given point 
$p=(x_0:y_0:z_0:w_0)\in \PP^3$ is in the surface $\Sc$, we just have to
substitute $x,y,z,w$ respectively by $x_0,y_0,z_0,w_0$ in ${\tt Z_1}$ and
check its rank; $p$ is on $\Sc$ if and only if the rank of
${\tt Z_1}$ does not drop. Of course some numerical aspects have to be taken
into account here, but the use of numerical linear algebra seems to be very
promising in this direction. Note also that 
such matrices, whose computation is very fast, are much more compact representations of implicit
equations compared to expanded polynomials which can have a lot of
monomials.

Even if the use of ${\tt Z_1}$ seems to be, in the opinion of the
authors, the best way to work quickly with implicit equations, the
previous algorithm also returns an explicit description of $P^\beta$, where $P$ is as usual an implicit
equation of $\Sc$ (always if $\Pc$ is locally a complete intersection). 
With the convention $\det(\Delta_i)=1$ if $\Delta_i$ 
does not exists, for $i=2,3$, $P^\beta$ can be computed as
the quotient $\frac{\det(\Delta_1)\det(\Delta_3)}{\det(\Delta_2)}$. Since $\det(\Delta_1)$ is
always a multiple of  $P^\beta$, we have
$\det(\Delta_1)=P^\beta Q$, where $Q$ is an homogeneous polynomial in
$\KK[\ux]$. It can be useful to notice that this extraneous factor $Q$
divides $\det(\Delta_2)$; in fact $\det(\Delta_2)=Q\det(\Delta_3)$.   

\subsection{Examples}
We illustrate our algorithm with four particular examples. It has been implemented in the software MAGMA which offers very
powerful tools to deal with linear algebra, and appears to be very efficient.     

\subsubsection{An example without base points}\label{nobp} 
 Consider the following example~:
 $$\left\{\begin{array}{ccl}
f_0 &= & s^2t, \\
  f_1 &= & t^2u,\\
  f_2 &=& su^2,\\
  f_3 & = & s^3+t^3+u^3.
\end{array}\right.$$
It is a parameterization without base points of a surface of degree
9. Applying our algorithm we find, in degree $\nu=2.3-2=4$ a matrix
${\tt Z_1}$ of size $15\times 24$ which represents our surface. 
Continuing the algorithm we finally obtain a matrix $\Delta_1$ of size
$15\times 15$, a matrix $\Delta_2$ of size  $9\times 9$ and a matrix $\Delta_3$
of size  $3\times 3$. Computing the quotient 
$\frac{\det(\Delta_1)\det(\Delta_3)}{\det(\Delta_2)}$ we obtain an
expanded implicit equation~:
$$\begin{array}{ll}
x^6z^3 + 3x^5y^2z^2 + 3x^4y^4z + 3x^4yz^4 + x^3y^6 +
        6x^3y^3z^3 + 3x^2y^5z^2 + 3x^2y^2z^5 - \\
        x^2y^2z^2w^3 +   3xy^4z^4 + y^3z^6.
\end{array}$$
Trying our algorithm empirically with $\nu=3$ we obtain a matrix $\Delta_1$ which
is not square, it has only 9 columns (instead of the waited 10). The
algorithm hence fails here, showing that in this case the degree
bound $2d-2$ is the lowest possible. 

\subsubsection{An example where the ideal of base points is saturated}
This example is taken from \cite{SeCh95} (see also
\cite{BCD02}), it is the  
parameterization of a cubic surface with 6 local complete intersection
base points~:
$$\left\{\begin{array}{ccl}
f_0 &= & s^{2} {t}+2 t^{3}+s^{2} {u}+4 {s} {t} {u}+4 t^{2} {u}+3 {s}
  u^{2}+2 {t} u^{2}+2 u^{3}, \\
  f_1 &= & -s^{3}-2 {s} t^{2}-2 s^{2} {u}-{s} {t} {u}+{s} u^{2}-2 {t}
  u^{2}+2 u^{3},\\
  f_2 &=& -s^{3}-2 s^{2} {t}-3 {s} t^{2}-3 s^{2} {u}-3 {s} {t} {u}+2 t^{2}
  {u}-2 {s} u^{2}-2 {t} u^{2},\\
  f_3 & = & s^{3}+s^{2} {t}+t^{3}+s^{2} {u}+t^{2} {u}-{s} u^{2}-{t}
  u^{2}-u^{3}.\end{array}\right.$$

The ideal $I=(f_0,f_1,f_2,f_3)$ is here saturated, so that we have $\nu_0=d-2=1$, and we hence can apply the algorithm with $\nu=1$. The 
 matrix ${\tt Z_1}$ is then square and the algorithm stops in step
2. It is given by~:
$$\left(\begin{array}{ccc}
                x       &     -z - w      &       y + w\\
                y & x - 2y + z - 2w       &    2y - z \\
                z  &        -x - 2w     &      y + 2w
                \end{array}\right).$$

\subsubsection{An example where the method of moving quadrics fails}
We here consider the example 3.2 in \cite{BCD02}. This example was
introduced to show how the method of moving quadrics (introduced in
\cite{SeCh95}) generalized in this paper to the presence of base
points may fail. Consider the parameterization
$$\left\{\begin{array}{ccl}
f_0 &= & su^2, \\
  f_1 &= & t^2(s+u),\\
  f_2 &=& st(s+u),\\
  f_3 & = & tu(s+u).
\end{array}\right.$$
We can directly apply our algorithm in degree $\nu=2.3-2=4$; 
we obtain a matrix
${\tt Z_1}$ of size $15\times 30$, matrices
$\Delta_1$ and $\Delta_2$ of size 
$15\times 15$, and a matrix $\Delta_3$ of size  $3\times 3$. An
expanded implicit equation is then $xyz + xyw - zw^2$. 

If now we take into account that the ideal $I=(f_0,f_1,f_2,f_3)$ has
base points, we can apply the algorithm with $\nu=3$. Even better, the
saturated ideal of $I$, denoted $I_\Pc$, is generated by $t(s+u)$
and $su^2$, showing that $\indeg(I_\Pc)=2$ (since $t(s+u)$ is
of degree 2), and hence that we can apply the algorithm with
$\nu=2$. In this case $\Delta_1$ is $6\times 6$, $\Delta_2$ is 
$3\times 3$ and the algorithm stops in step $(e)$.

Remark that to have the best bound possible for the integer $\nu$ we
have to compute the initial degree of $I_\Pc$, which is not always
obvious. However, it should be possible to determine some classes of
parameterizations which satisfy  certain geometric properties and for
which we know in advance the initial degree of the saturation of $I$.

\subsubsection{An example with a fat base point}

With this example we would like to illustrate the theorem 
\ref{theorem surface} when there is a base point minimally generated by three polynomials. Consider the parameterization given by
$$\left\{\begin{array}{ccl}
f_0 &= &s^3-6s^2t-5st^2-4s^2u+4stu-3t^2u , \\
  f_1 &= &-s^3-2s^2t-st^2-5s^2u-3stu-6t^2u ,\\
  f_2 &=&-4s^3-2s^2t+4st^2-6t^3+6s^2u-6stu-2t^2u ,\\
  f_3 & = & 2s^3-6s^2t+3st^2-6t^3-3s^2u-4stu+2t^2u  .
\end{array}\right.$$

The ideal $I=(f_0,f_1,f_2,f_3)$ defines exactly one \emph{fat} base point $p$ which is defined by $(s,t)^2$. Therefore any implicit equation of our parameterized surface is of degree $d^2-e_p=9-4=5$, and the degree of the determinant of the complex ${(\Zc_\bullet)}_{[\nu_0]}$ is of degree $d^2-d_p=9-3=6$. Applying our algorithm with $\nu_0=2(d-1)-2=2$ we find that the matrix ${\tt Z_1}$ is already square, of size $6\times 6$. Expanding its determinant one obtains a product of a degree 5 irreducible polynomial (an implicit equation of our surface) and a linear (since $e_p-d_p=1$) irreducible polynomial.

\section{Implicitization of rational parametric hypersurfaces}

We know turn to the general problem of hypersurface implicitization,
always using approximation complexes. In this section we prove new
results to solve this problem under suitable conditions, and obtain
the proof of theorem \ref{theorem surface} as a particular case.  

We are going to consider hypersurfaces obtained as the closed image of
maps from $\PP^{n-1}$ to $\PP^n$, where $n\geq 3$. Consequently we set
$A:=\KK[X_1,\ldots,X_n]$, which will be the ring of the polynomials of
the parameterizations we will consider. We also denote $\im 
:=(X_{1},\ldots ,X_{n})$, and set $\hbox{---}^{\vee}:=\hom_{A}
(\hbox{---},A[-n])$ and $\hbox{---}^{*}:=\homgr_{A}(\hbox{---},A/\im
)$. Finally $\omega_{-}$ will denote the associated canonical module
(see \cite{BrHe93}). Our first result is a technical lemma on the
cycles of certain Koszul complexes.   

\begin{lemma}\label{Zcohom} Let $f_{0},\ldots ,f_{n}$ be $n+1$ polynomials of positive
  degrees $d_{0},\ldots ,d_{n}$, $I$ the 
  ideal generated by them, ${\bf   K}_{\bullet}(f;A)$ be the Koszul
  complex of the $f_{i}$'s on $A$ and denote by $Z_{i}$ and $H_{i}$
  the $i$-th cycles and the $i$-th  homology modules of  ${\bf
  K}_{\bullet}(f;A)$, respectively. Denoting $\s := d_{0}+\cdots
  +d_{n}$, if $\dim (A/I)\leq 1$  then we have 
\begin{itemize}
\item $H_{i}\not= 0$ for $i=0,1$, $H_{i}=0$ for $i>2$ and $H_{2}$ is
    zero if and only if $\dim (A/I)=0$. If $\dim (A/I)=1$,
    $H_{2}\simeq \om_{A/I}[n-\s ]$. 
\item If  $n\geq 3$,  then 
$$ H^{i}_{\im}(Z_{p})\simeq\left\{\begin{array}{cl}
 0 & \hbox{for}\ i=0,1  \\
 H^{0}_{\im}(H_{i-p})^{*}[n-\s ] & \hbox{for}\ i=2 \\
 H_{i-p}^{*}[n-\s ] & \hbox{for}\ 2<i<n \\
 Z_{n-p}^{*}[n-\s ] & \hbox{for}\ i=n.
\end{array}\right.$$
\end{itemize}
\end{lemma}

\begin{proof} The first point is classical, see e.g. \cite{BrHe93}
  1.6.16 and 1.2.4. For the second point we consider the truncated
  complexes: 
$$
{\bf K}_{\bullet}^{>p}:\quad \quad 0\rightarrow K_{n+1}\rightarrow \cdots
\rightarrow K_{p+1}\rightarrow Z_{p}\rightarrow 0.$$
If $p> 2$ this complex is exact, and gives rise to a 
spectral sequence $H^{\bullet}_{\im}({\bf
  K}_{\bullet}^{>p})\Rightarrow 0$, which is at level 1~: 
\begin{diagram}
0 & & \cdots & & 0 & & H^{0}_{\im}(Z_{p})\\
\vdots& & &  & \vdots & & H^{1}_{\im}(Z_{p})\\
\vdots&& & & \vdots& & \vdots\\
0&&\cdots&&0&&H^{n-1}_{\im}(Z_{p})\\
H^{n}_{\im}(K_{n+1})& \rTo^{
  H^{n}_{\im}(\partial_{n+1})}&\cdots&\rTo^{H^{n}_{\im}(\partial_{p+2})}
&H^{n}_{\im}(K_{p+1})&\rTo^{H^{n}_{\im}(\partial_{p+1})}&H^{n}_{\im}(Z_{p}),\\ 
\end{diagram}
and modulo the identifications
$$H^{n}_{\im}(K_{\bullet})\simeq (K_{\bullet}^{\vee})^{*}\simeq
(K^{\bullet}[-n])^{*}\simeq (K_{n+1-\bullet}[\s -n])^{*},
$$
the last line becomes
$$K_{0}^{*}[n-\s ]\xrightarrow{\partial_{1}^{*}} K_{1}^{*}[n-\s]
\xrightarrow{\partial_{2}^{*}} 
\cdots \xrightarrow{\partial_{n-p}^{*}} K_{n-p}^{*}[n-\s ] 
\rightarrow H^{n}_{\im}(Z_{p}).$$
It follows 
\begin{itemize}
\item $H^{i}_{\im}(Z_{p})=0$ for $i<p$ and for $p+2<i<n$,
\item $d^{>p}_{n-p+1}: H_{0}^{*}[n-\s] \xrightarrow{\simeq} 
H^{p}_{\im}(Z_{p})$,
\item $d^{>p}_{n-p}: H_{1}^{*}[n-\s]\xrightarrow{\simeq}
 H^{p+1}_{\im}(Z_{p})$,
\item $d^{>p}_{n-p-1}: H_{2}^{*}[n-\s] \xrightarrow{\simeq}
 H^{p+2}_{\im}(Z_{p})$,
\end{itemize}
where the upper index indicates from which spectral sequence the
isomorphism comes, and the lower index indicates at which level of the
spectral sequence the map is obtained. We also obtain a short exact sequence
$$K_{n-p-1}^{*}[n-\s ]\xrightarrow{\partial_{n-p}^{*}}
K_{n-p}^{*}[n-\s ]\rightarrow H^{n}_{\im}(Z_{p})\rightarrow 0$$
that gives the isomorphism $H^{n}_{\im}(Z_{p})\simeq Z_{n-p}^{*}[n-\s ]$.

For $p=2$, we have a spectral sequence $H^{\bullet}_{\im}({\bf
  K}_{\bullet}^{>2})\Rightarrow H^{1}_{\im}(H_{2})$ which shows as in
  the previous cases that $
H^{n}_{\im}(Z_{2})\simeq Z_{n-2}^{*}[n-\s ]$, 
and that 
\begin{itemize}
\item $H^{i}_{\im}(Z_{2})=0$ for $4<i<n$,
\item $d^{>2}_{n-2}:H_{1}^{*}[n-\s] \xrightarrow{\simeq}
 H^{3}_{\im}(Z_{2})$,
\item $d^{>2}_{n-3}:H_{2}^{*}[n-\s]\xrightarrow{\simeq}
H^{4}_{\im}(Z_{2})$.
\end{itemize}
It also provides an exact sequence 
\begin{equation}\label{star}
0\rightarrow H^{1}_{\im}(Z_{2})\xrightarrow{can}
H^{1}_{\im}(H_{2})\xrightarrow{\tau}
H_{0}^{*}[n-\s ]\xrightarrow{d^{>2}_{n-1}}
H^{2}_{\im}(Z_{2})\rightarrow 0,
\end{equation}
where $\tau$ is the transgression map of the spectral sequence. If
$H_{2}\not= 0$, we also look at the complex
$$
{\bf K}_{\bullet}^{\leq 2}:\quad 0\rightarrow Z_{2}\rightarrow K_{2}
\rightarrow K_{1}\rightarrow I\rightarrow 0,$$
whose only homology is $H_{1}$ and the corresponding spectral sequence
$H^{\bullet}_{\im}({\bf 
  K}_{\bullet}^{\leq 2})\Rightarrow H^{\bullet}_{\im}(H_{1})$. 
Noticing that
$H^{i}_{\im}(H_{1})=0$ for $i>1$, we get that
$H^{0}_{\im}(Z_{2})=H^{1}_{\im}(Z_{2})=0$ (this is also easily
obtained by splitting ${\bf K}_{\bullet}^{\leq 2}$ into two short
exact sequences, and taking cohomology). Together with \eqref{star} and
local duality, which gives $H^{1}_{\im}(H_{2})\simeq \om_{H_{2}}^{*}\simeq
(A/I^{{\rm sat}})^{*}[n-\s ]$ where $I^{{\rm sat}}$ denotes the saturation of the ideal $I$, we get the asserted
isomorphism for $H^{2}_{\im}(Z_{2})$.

Note that we also get an isomorphism $H^{0}_{\im}(H_{0})^{*}[n-\s
]\simeq H^{2}_{\im}(Z_{2})\simeq
H^{0}_{\im}(H_{1})$, and an exact sequence: 
$$
0\rightarrow H^{1}_{\im}(H_{1})\rightarrow H^{3}_{\im}(Z_{2})
\rightarrow H^{1}_{\im}(I)\rightarrow 0.$$

Finally for $p=1$, the exact sequence
$$
{\bf K}_{\bullet}^{\leq 1}:\quad 0\rightarrow Z_{1}\rightarrow K_{1}
\rightarrow I\rightarrow 0$$
shows that $H^{i-2}_{\im}(H_{0})\simeq H^{i-1}_{\im}(I)\simeq
H^{i}_{\im}(Z_{1})$ for $i<n$, so that
$H^{0}_{\im}(Z_{1})=H^{1}_{\im}(Z_{1})=0$, $H^{2}_{\im}(Z_{1})\simeq
H^{0}_{\im}(H_{0})\simeq H^{0}_{\im}(H_{1})^{*}[n-\s
]$, $H^{3}_{\im}(Z_{1})\simeq
H^{1}_{\im}(H_{0})\simeq H_{2}^{*}[n-\s ]$ and $H^{i}_{\im}(Z_{1})=0$
for $3<i<n$.  

The spectral sequence $H^{\bullet}_{\im}({\bf
  K}_{\bullet}^{>1})\Rightarrow H^{\bullet}_{\im}(H_{\bullet})$
  (because only $H^{1}_{\im}(H_{2})$,  $H^{0}_{\im}(H_{1})$ and
  $H^{1}_{\im}(H_{1})$ may not be zero) gives
  $H^{n}_{\im}(Z_{1})\simeq Z_{n-1}^{*}[n-\s ]$, and this concludes
  the proof. 
\end{proof}

\begin{remark}
Note that the last spectral sequence of the proof also identifies
$H^{2}_{\im}(Z_{1})$ in another way by providing the exact sequence
$$
0\rightarrow H^{1}_{\im}(H_{1})\xrightarrow{\tau '} H_{1}^{*}[n-\s ]
\xrightarrow{d_{n-1}^{>1}} H^{2}_{\im}(Z_{1})\rightarrow 0
$$
which shows that if $\dim (A/I)=1$
$$
H_{1}/H^{0}_{\im}(H_{1})\simeq\om_{H_{1}}[n-\s ].
$$
Thus in this case $H_{1}/H^{0}_{\im}(H_{1})$ has a symmetric free
resolution. Also, the spectral sequences derived from the complexes
${\bf K}_{\bullet}^{\leq p}$ provide exact sequences
$$
0\rightarrow H^{0}_{\im}(H_{1})\rightarrow H^{p}_{\im}(Z_{p})\rightarrow 
H^{1}_{\im}(H_{2})\rightarrow 0
$$
for $p<n$, and 
$$
0\rightarrow H^{1}_{\im}(H_{1})\rightarrow H^{p+1}_{\im}(Z_{p})
\rightarrow
H^{0}_{\im}(H_{0})\rightarrow 0$$
for  $p<n-1$, as well as an isomorphism $H^{p+2}_{\im}(Z_{p})\simeq
H^{1}_{\im}(H_{0})$ for $p<n-2$. 
\end{remark}

As we have already done in previous sections, we now consider
approximation complexes. To do this, we introduce new variables
$T_{1},\ldots ,T_{n+1}$ which represent the homogeneous coordinates of
the target $\PP^n_\KK$ of a given parameterization. Let
$f_0,\ldots,f_n$ be $n+1$ homogeneous polynomials in $A$ of the same
degree $d\geq 1$. Denoting by $Z_i$ the $i^{\mathrm{th}}$-cycles of
the Koszul complex of the $f_i$'s on $A$, we set
$\Zc_i:=Z_i[id]\otimes_A A[\underline{T}]$, where $[\hbox{---}]$
stands for the degree shift in the $X_{i}$'s and $(\hbox{---})$ for
the one in the $T_{i}$'s. It appears that the differentials  
$v_\bullet$ of the Koszul complex $K_{\bullet}(T_{1},\ldots
,T_{n+1};A[\underline{T}])$ induce maps between the $\Zc_i$'s, and hence we can
define the approximation complex (note that $\Zc_{n+1}=0$) 
$$(\Zc_\bullet,v_\bullet): 0\rightarrow \Zc_{n}(-n) \xrightarrow{v_n}
\ldots \xrightarrow{v_3} \Zc_2(-2) \xrightarrow{v_2} \Zc_1(-1)
\xrightarrow{v_1} \Zc_0=A[\underline{T}].$$ 
It is a naturally a bi-graded complex, and it is easy to check that 
$H_{0}(\Zc_{\bullet})\simeq \sym_{A}(I)$. We now give some acyclicity
criterions in case $I=(f_0,\ldots,f_n)$ define isolated points in
$\Proj(A)$. Recall that a sequence $x_1,\ldots,x_n$ of elements in a
ring $R$ is said to be a \emph{proper sequence} if
$$x_{i+1}H_{j}(x_1,\ldots,x_i;A)=0 \ \ \mathrm{for} \ i=0,\ldots,n-1 \
\mathrm{and} \ j>0,$$ where the $H_j$'s denote the homology groups of
the associated Koszul complex. 

\begin{lemma}\label{Zacyclicity} Suppose that $I=(f_{0},\ldots
  ,f_{n})$ is of codimension at least $n-1$, $\KK$ is infinite, and let $\Pc
  :=\Proj (A/I)$. Then the following are equivalent~: 
\begin{enumerate}
\item[\rm (1)] ${\Z}_{\bullet}$ is acyclic,
\item[\rm (2)] ${\Z}_{\bullet}$ is acyclic outside $V(\im )$,
\item[\rm (3)] $I$ is generated by a proper sequence,
\item[\rm (4)] $\Pc$ is locally defined by a proper sequence,
\item[\rm (5)] $\Pc$ is locally defined by $n$ equations.
\end{enumerate}
\end{lemma}
\begin{proof} By \cite{HSV83L}, (1)$\Leftrightarrow$(3) and
(2)$\Leftrightarrow$(4). Moreover (1)$\Rightarrow$(2) and
(3)$\Rightarrow$(4). We will show that
(4)$\Rightarrow$(5)$\Rightarrow$(3). 

Assume (4). For each $p\in\Pc$ there exists a non empty open set $\Omega_{p}$ in
$\KK^{(n+1)^{2}}$ such that if $(a_{ij})\in \Omega_{p}$ and
$g_{i}:=\sum_{j}a_{ij}f_{j}$, the $g_{i}$'s form a proper sequence
locally at $p$. Taking $(a_{ij})\in\Omega:=\cap_{p\in\Pc} \Omega_{p}$, this gives
$g_{1},\ldots ,g_{n+1}$ so that the $g_{i}$'s form a proper
sequence outside $V(\im )$. We may also assume, by shrinking $\Omega$
if necessary, that $g_{1},\ldots ,g_{n-1}$ is a regular sequence, so that
$g_{1},\ldots ,g_{n}$  is  clearly a proper sequence. 
Now $g_{1},\ldots ,g_{n+1}$ is proper if and only if
$g_{n+1}$ annihilates $H_{1}(g_{1},\ldots ,g_{n};A)\simeq
\om_{A/J}$ (locally outside $V(\im )$, {\it a priori}), where $J$ denotes the saturated ideal of $(g_1,\ldots,g_n)$ w.r.t. $\im$. But $\ann_{A}(\om_{A/J})=J$, so that $g_{n+1}\in J$, which proves (5).      

Now assume (5). One considers (in the same spirit as above) $g_{1},\ldots
,g_{n+1}$ so that $g_{1},\ldots ,g_{n-1}$ is a regular sequence and  
$g_{1},\ldots ,g_{n}$ defines $\Pc$. Then $g_{n+1}\in I_{\Pc}$,
and therefore annihilates $\om_{A/I_{\Pc}}$, so that the $g_{i}$'s
forms a proper sequence.
\end{proof}

We are now able to state our main result on the hypersurface implicitization problem.

\begin{theorem}\label{TH} Let $n\geq 3$. Let $I$ be the ideal
  $(f_{0},\ldots ,f_{n})$, $\Pc :=\Proj (A/I)$ and $I_{\Pc}$ the 
saturation of $I$ w.r.t. $\im$. Assume that $\KK$ is infinite and that $\dim \Pc \leq 0$, that is
$\Pc$ define a finite number of base points in $\Proj(A)$, possibly
empty. Then we have  
\begin{itemize}
\item ${\Z}_{\bullet}$ is acyclic if and only if $\Pc$ is
locally defined by $n$ equations.
\item  Let $\nu \geq \nu_{0}:=(n-1)(d-1)-\indeg
(I_{\Pc})$. Then $H^{0}_{\im}(\sym_{A}(I))_{[\nu]}=0$. Moreover $({\cal
  Z}_{\bullet})_{[\nu]}\otimes_{\KK [\underline{T}]}\KK
(\underline{T})$ is acyclic if and only if  ${\Z}_{\bullet}$ is
acyclic. 

\item Let $\nu \geq \nu_{0}$ and assume that $\Pc$ is locally defined by $n$
equations. Then $D:=\det (({\cal   Z}_{\bullet})_{[\nu]})$ is a non zero
homogeneous element of $\KK[\underline{T}]$, independent of $\nu$ (modulo
$\KK^{\times}$), of degree $d^{n-1}-\sum_{p\in\Pc}d_p$. Denoting by
$\mathcal{H}$ the closed image of the rational map $\phi :{\bf
  P}^{n-1} \rightarrow {\bf P}^{n}$, $D=H^{\deg(\phi)}G $  where
$H$ is an implicit equation of $\mathcal{H}$. Moreover $G\in \KK^\times$ if and only if
$\Pc$ is locally of linear type, and if and only if $\Pc$ is locally a
complete intersection. 
\end{itemize}
\end{theorem}

Before giving the proof of the theorem we recall that, by definition,
$\Pc$ is said to be locally of linear type if
$\Proj(\sym_A(I))=\Proj(\mathrm{Rees}_A(I))$. By  \cite{BuJo02}
theorem 2.5, $\deg(H^{\deg(\phi)})=d^{n-1}-\sum_{p\in\Pc}e_p$, and
consequently we always have
$\deg(G)=\sum_{p\in\Pc}(e_p-d_p)$. Also  $e_p\geq d_p$ with equality
if and only if $\Pc$ is locally a complete intersection at $p\in\Pc$. 

\begin{proof} The first point follows from lemma \ref{Zacyclicity}. For the second point we consider the
two spectral sequences associated to the 
double complex $H^{\bullet}_{\im}(\Z_{\bullet})$, both abouting to the
hypercohomology of $\Z_{\bullet}$. One of them abouts at level two
with:

$$
{_{2}{'}E}^{p}_{q}={_{\infty}{'E}}^{p}_{q}=\left\{ 
\begin{array}{cl}
H^{p}_{\im}(H_{q}(\Z_{\bullet})) &\hbox{for}\ p=0,1\
  \hbox{and}\ q>0 \\
  H^{p}_{\im}(\sym_{A}(I)) &\hbox{for}\ q=0\\
0 &\hbox{else}.\\
\end{array}\right.
$$
The other one gives at level one:
$$
{_{1}{''E}}^{p}_{q}=H^{p}_{\im}(Z_{q})[qd]\otimes_{A}A[\underline{T}](-q).
$$

By lemma \ref{Zcohom}, $H^{p}_{\im}(Z_{q})=0$ for $p<q$ and for $p<2$. This shows
that $H^{p}_{\im}(H_{q}(\Z_{\bullet}))=0$ for $q>0$
except possibly for $p=q=1$, by comparing the two spectral
sequences. Therefore, $H_{q}(\Z_{\bullet})=0$ for $q>1$ and 
$H^{0}_{\im}(H_{1}(\Z_{\bullet}))=0$. 

If $p>2$
$$
{_{1}{''E}}^{p}_{p}\simeq H^{*}_{0}[n-(n+1-p)d]\otimes_{A}A[\underline{T}](-p)
$$
so that $({_{1}{''E}}^{p}_{p})_{\nu}=0$ if $\nu >(n-2)d-n$. 
Also, 
$$
({_{1}{''E}}^{2}_{2})\simeq (I_{\Pc}/I)^{*}[n-(n-1)d]\otimes_{A}A[\underline{T}](-2),
$$
so that $({_{1}{''E}}^{2}_{2})_{\nu}=0$ for $\nu
>(n-1)d-n-\indeg (I_{\Pc}/I)$. 
\medskip

Therefore, if $\nu \geq \nu_{0}$, $({_{1}{''E}}^{p}_{q})_{\nu}=0$ for
$p\leq q$. Note also that we have the equalities 
$$\min \{d,\indeg (I_{\Pc}/I)\}=\min \{ d,\indeg (I_{\Pc})\}=\indeg (I_{\Pc}).$$ By comparing with the other spectral
sequence, we have $H^{1}_{\im}(H_{1}(\Z_{\bullet}))_{\nu}=0$ and
$H^{0}_{\im}(\sym_{A}(I))_{\nu}=0$.  

As $H^{0}_{\im}(\sym_{A}(I))_{\nu}=0$ for $\nu \geq \nu_{0}$ we have
$$
\ann_{\KK [\underline{T}]}(\sym_{A}(I)_{\nu})=\ann_{\KK
  [\underline{T}]}(\sym_{A}(I)_{\nu_{0}})
$$
for any $\nu\geq \nu_{0}$ (see for instance \cite{BuJo02}, the proof
of 5.1), so that this module is torsion if and only if $\Pc$ is 
locally defined by at most $n$ equations (because $\ann_{\KK
  [\underline{T}]}(\sym_{A}(I)_{\nu})$ is torsion for $\nu\gg 0$ if
and only if $I$ is defined by $<n+1$ equations outside $V(\im
)$; one may also use the study of the minimal primes of $\sym_{A}(I)$
in \cite{HuRo86}). Also in the case where $\Pc$ is locally defined by
at most $n$ equations, the divisor associated to this module is
independant of $\nu\geq \nu_{0}$. This finishes the proof of point
two.    

We come to the third point. Note that we have just  shown that $D$ is
independant of $\nu$ (up to an element of $\KK^{\times}$). To compute
the degree of $D$, we may then take  $\nu \gg 0$.
 
The matrices of the maps of $\Zc_{\bullet}$ have entries
which are linear forms in the $T_{i}$'s, so that the determinant of
${(\Zc_\bullet)}_{[\nu]}$ is a form in the $T_{i}$'s of
degree  $$\delta :=\sum_{i=1}^n
(-1)^{i+1}i\dim_{\KK}({Z_i}_{[\nu+id]}).$$

We have canonical exact sequences, $i=0,\ldots,n$,
$$0\rightarrow Z_{i+1}\rightarrow K_{i+1}
\rightarrow B_{i} \rightarrow 0$$
and 
$$0\rightarrow B_{i}\rightarrow Z_{i}
\rightarrow H_{i} \rightarrow 0$$
which shows that $\delta$ can be expressed in terms
of the Hilbert polynomials of the $K_{i}$'s and the $H_{i}$'s.
The contribution of the $K_{i}$'s only depends on $n$ and $d$, and the
one of the $H_{i}$'s only comes from $Z_{1}$ and $Z_{2}$ and is : 
$$
(H_{0})_{\nu +d}-2((H_{1})_{\nu +2d}-(H_{0})_{\nu +2d})
$$
as $\deg H_{1}=2\deg \Pc$ (one may use for example that
$(H_{0})_{\nu}-(H_{1})_{\nu}+(H_{2})_{\nu}=0$ for $\nu\gg 0$ and that
$H_{2}\simeq \om_{R/I_{\Pc}}$ up to a degree shift), this contribution
is equal to $-\deg \Pc$ for $\nu \gg 0$. The contribution of the
$K_{i}'s$ is $d^{n-1}$ (for $\nu \gg 0$) in the case where the
$H_{i}$'s are $0$ for $\nu \gg 0$, and we get $\delta =d^{n-1}-\deg
\Pc$. 

Let $ X:=\proj ({\rm Rees}_{I}(A))\subseteq Y:=\proj
  (\sym_{A}(I))\subset {\bf P}^{n-1}\times {\bf P}^{n}$. The scheme
  $X$ is the closure of the graph of $\phi$ in ${\bf P}^{n-1}\times {\bf
  P}^{n}$. Therefore, the closed image of $\phi$ is $p_{2}(X)$,
  and  $p_{2}(X)\subseteq p_{2}(Y)$. As $D\not= 0$,
  $p_{2}(Y)\not= {\bf P}^{n}$. If $\Pc$ is locally of linear type,
  then $X=Y$, so that $G\in \KK^{\times}$ which implies that $d_p=e_p$
  for any $p\in \Pc$ and therefore $\Pc$ is locally a 
  complete intersection. If $\Pc$ is not locally a complete
  intersection at $p\in \Pc$, $d_p\not= e_p$ so that $\deg G>0$ and
 $X\not= Y$.
 Note that these last facts also follows from the study of the minimal
  primes associated to $\sym_{A}(I)$ in  \cite{HuRo86}.
\end{proof}

This theorem yields easily an algorithm to implicitize parameterized hypersurfaces under suitable assumptions, using only linear algebra routines. The case $n=3$ have been completely explicited in section \ref{surf}.

\bibliographystyle{plainnat}
\bibliography{volos2b.bbl}

\end{document}